\documentclass[11pt]{amsart} 

\usepackage{amsmath, amssymb} 
\usepackage{amscd}            
\usepackage{amsxtra}          
\usepackage{upref}            

\usepackage[mathscr]{eucal}

\theoremstyle{plain}

\newtheorem{thm}{Theorem}
\newtheorem{proposition}[thm]{Proposition}

\newtheorem{corollary}[thm]{Corollary}

\theoremstyle{remark}
\newtheorem{remark}[thm]{Remark}

\numberwithin{equation}{section}
\numberwithin{thm}{section}
\newcommand{\be}%
  {\protect\setcounter{equation}{\value{subsubsection}}}  
\newcommand{\ee}%
  {\protect\setcounter{subsubsection}{\value{equation}}}

\newcommand{\bk}{{\bar k}}
\newcommand{\bX}{X_{\bk}}
\newcommand{\bA}{A_{\bk}}
\newcommand{\rat}{{\rm rat}}
\newcommand{\alg}{{\rm alg}}
\renewcommand{\hom}{{\rm hom}}
\newcommand{\ltor}{{\ell{\rm -tors}}}

\newcommand{\Gr}{\operatorname{Gr}}

\newcommand{\isom}{\cong}          

\newcommand{\CH}{\operatorname{CH}}

\renewcommand{\P}{{\mathbb{P}}}          
\newcommand{\Q}{{\mathbb{Q}}}          

\newcommand\simto{\overset{\sim}{\longrightarrow}}
\newcommand{\tensor}{\otimes}
\newcommand{\Z}{{\mathbb Z}}

\begin{document}


\title {Griffiths groups of \\ supersingular Abelian varieties}
\author{B.~Brent Gordon}  
\address{Department of Mathematics \\ 
         University of Oklahoma \\
         Norman,\enspace OK  73019 \\ USA}
\email{bgordon@math.ou.edu}

\author{Kirti Joshi}
\address{Department of Mathematics \\ 
         University of Arizona \\
         Tucson,\enspace AZ  85721 \\ USA}
\curraddr{School of Mathematics \\
          Tata Institute for Fundamental Research \\
          Homi Baba Road \\
          Mumbai 400005, INDIA}
\email{kirti@math.tifr.res.in}



\maketitle


\markboth{GORDON, JOSHI}{GRIFFITHS GROUPS}


\section{Introduction}
\label{intro}

Let $k$ be a finite field of characteristic $p>0$. We fix an algebraic 
closure $\bk$ of $k$. For any finite extension $k'/k$ we write $G_{k'}$
for the Galois group of $\bk/k'$. Let $X/k$ be a smooth,  projective
variety over $k$. We will write $\bX$ for $X\times_k \bk$. Let
$Z^r(\bX)$ be the group of codimension $r\geq 0$ cycles on $\bX$. Let
$Z^r(\bX)_\rat$, $Z^r(\bX)_\alg$ and $Z^r(\bX)_\hom$ be the  subgroups
of codimension~$r$ cycles which are rationally, respectively
algebraically, respectively homologically equivalent to zero on~$\bX$. 
We will write $\CH^r(\bX)_\alg\subset \CH^r(\bX)_\hom\subset \CH^r(\bX)$
for the corresponding groups modulo the subgroup of cycles rationally
equivalent to zero.  Then the quotient $\CH^r(\bX)/\CH^r(\bX)_\hom$ is
finitely generated modulo torsion, and the Tate conjecture predicts the
rank of this group to be the order of vanishing of a suitable
$L$-function \cite{tate65}.  

	The quotient $\CH^r(\bX)_\hom/\CH^r(\bX)_\alg$ was first
investigated by P.~Griffiths and is called the {\em Griffiths group\/}
of~$X_\bk$. In every example where the structure of this group is
known, when it is not trivial it is quite subtle. We will write
$\Gr^r(\bX)$ for the group $\CH^r(\bX)_\hom/\CH^r(\bX)_\alg$, and refer
to it as the Griffiths group of codimension $r$ cycles.  Recently, Chad
Schoen has investigated the structure of the Griffiths group of
varieties over the algebraic closure of finite fields (see
\cite{schoen95}, \cite{schoen99}) and has shown that these groups can be
infinite in several interesting situations.  This note is inspired by
these papers of Schoen.

Recall that over a algebraically closed field of characteristic~$p$ a
{\em supersingular abelian variety\/} may be characterized by being
isogenous to a product of supersingular elliptic curves (see
\cite[Theorem 4.2]{oort74}), where an elliptic curve is said to be
supersingular if it possesses no geometric points of order
exactly~$p$.  The purpose of this note is to prove that the Griffiths
group of any supersingular abelian variety is at most a $p$-primary
torsion group (see Theorem~\ref{main}).  In \cite[Theorem~14.4, page
45]{schoen95}, Schoen had proved this assertion for the triple product
of the Fermat cubic and $p\equiv 2\pmod 3$.  Our result also applies
to all supersingular Fermat threefolds (see Theorem~\ref{fermat}).
One of the key ingredients in our proof is the work of N.~Fakhruddin
(see \cite{fakhruddin95}).  We hope to study the $p$-primary torsion
in a forthcoming work.
 
We would like to thank C.~Schoen and N.~Fakhruddin for correspondence
and comments. We are also grateful to F.~Oort and the referee for a
very careful reading and helpful suggestions which have improved the
readability of this paper.

\section{Beauville's conjecture}\label{beauvillesection}
Suppose $A$ is an abelian variety over $k$. Then by results of Mukai,
Beauville, Deninger-Murre (see \cite{mukai81}, \cite{beauville86} and
\cite[Theorem 2.19, page 214]{deninger91}), the Chow groups of $A$
when tensored with $\Q$ admit a finite decomposition:
\begin{equation}\label{beauville-decomp}
\CH^i(\bA)\tensor\Q = \bigoplus_{j}\CH^i_j(\bA) ,
\end{equation}
 where $\CH^i_j(\bA)$ is the subset of $\CH^i(\bA)\tensor\Q$ on which
the flat pull-back of multiplication in~$\bA$ by an integer~$m$ acts
as multiplication by~$m^{2i-j}$. 

\begin{remark}
  To obtain the decomposition \eqref{beauville-decomp}, it is not
neccessary to admit $\Q$-coefficients, but it suffices to invert
integers which arise as denominators in the Riemann-Roch Theorem.  In
particular such a decomposition is available over a suitable
localization of~$\Z$.  From now on we will work with this refined
decompostion over a suitable localisation of $\Z$.
\end{remark}

	Using the work of Soul\'e (see \cite{soule84}), K\"unnemann
(see \cite[Theorem 7.1, page 99]{kunneman93}) proves that, except possibly
for $j=0$, all the remaining components of
\eqref{beauville-decomp} are torsion.

\begin{thm}\label{purity}
Suppose $A$ is an abelian variety over a finite field  $k$. Then for all
$i\geq 0$:
\begin{equation}
\CH^i(\bA)\tensor\Q = \CH^i_0(\bA)
\end{equation}
\end{thm}
This result is essentially a consequence of the fact that the motive
of an abelian variety is pure in the sense
of~\cite[Definition 3.1.1, page 331]{soule84}.  In particular, all the
groups $\CH^i_j(\bA)$ except possibly $\CH^i_0(\bA)$ must be torsion.

	Another version of this result, valid over any algebraically
closed field (of positive characteristic), was proved by N.~Fakhruddin
(see \cite{fakhruddin95}).

\begin{thm}\label{twostep}
  When $X$ is a supersingular abelian variety over an algebraically
closed field of characteristic~$p$, then $\CH^i_j(X) = 0$ for $j \neq
0,\, 1$.
\end{thm}

\section{A result of Fakhruddin}
 We will also need the following result of N.~Fakhruddin (see
\cite{fakhruddin95}). 
 
\begin{thm}\label{najmuddin}
Suppose $A$ is a supersingular abelian variety over an algebraically
closed field of characteristic~$p$. Then the  restriction of the
cycle class map to $\CH^d_0(A)$ induces an isomorphism:
\[ 
\CH^d_0(A)\tensor \Q_\ell \simto H^{2d}_{et}(A,\Q_\ell(d)) .
\]
\end{thm}

\begin{remark}
In particular, $\CH^i_1(X)$ is the kernel of the cycle class map when
$X$ is a supersingular abelian variety over any algebraically closed
field.
\end{remark}

	The idea behind the proof of this theorem is that the subgroup
$B^d(A)$ of $\CH^d_0(A)$ generated by classes of abelian subvarieties
of~$A$ of codimension~$d$ on the one hand coincides with $\CH^d_0(A)$
and on the other hand, after tensoring with~$\Q_\ell$ maps
isomorphically onto $H^{2d}_{et}(A,\Q_\ell(d))$.  See
\cite{fakhruddin95} for the details.

\section{Abel-Jacobi Maps}\label{ajmap}
 In \cite{bloch79}, Bloch constructed an Abel-Jacobi mapping
\[
\lambda_i: \CH^i(\bX)_\ltor \to
H^{2i-1}_{et}(\bX,\Q_\ell/\Z_\ell(i)) 
\]
 We will need the following results about this map.

\begin{thm}\label{merkurev-suslin}
{\rm (see \cite[Corollary4, page 775]{colliot-thelene83},
 \cite[Section 8]{merkurev83}) } Let $X/k$ be any smooth, projective
 variety and $\ell\neq p$.  Then the map $\lambda_2$ is injective.
\end{thm}

\begin{thm}\label{suwa}
{\rm (see \cite[Th\'eor\`eme 4.7.1, page 87]{suwa88}) } Let $A/k$ be a
 supersingular abelian variety, and $\ell\neq p$.  Then the
 restriction $\lambda'_i$ of $\lambda_i$ to
 $\CH^i(\bA)_{\alg,\,\ltor}$ is surjective, and bijective if $i=1,\,2$
 or $\dim A$.
\end{thm}

\section{Supersingular Abelian Varieties}

We are now in a position to prove the main theorem of this note.

\begin{thm}\label{main}
Let $\bA$ be a supersingular abelian variety over the algebraic closure
of a finite field. Then $\Gr^2(\bA)$ is at most a $p$-primary torsion
group.
\end{thm}

\begin{proof}
  Under the hypotheses of this theorem, the results in
section~\ref{beauvillesection} together with Theorem~\ref{najmuddin}
imply that
\[
 \CH^i_1(\bA) = \CH^i(\bA)_\hom 
\]
 is a torsion group.  Thus $\CH^i(\bA)_\alg$ and $\Gr^i(\bA)$ are also
torsion groups.  So to prove $\Gr^2(\bA)$ has no $\ell$-primary torsion
for $\ell\neq p$ it will suffice to prove that 
$
\Gr^2(\bA)\tensor\Z_\ell=0
$,
or, equivalently, that 
\begin{equation}\label{torgps}
\CH^2(\bA)_\hom\tensor\Z_\ell=\CH^2(\bA)_\alg\tensor\Z_\ell.
\end{equation}

   We use the Bloch-Abel-Jacobi mapping (section~\ref{ajmap})
to prove \eqref{torgps}.  Consider the commutative diagram:
\begin{equation}\label{diagram}
\begin{CD}
\CH^2(\bA)_\ltor  @>\lambda_2>>  H^3_{et}(\bA,\Q_\ell/\Z_\ell(2))\\
@AAA                                  @AA{\lambda'_2}A\\
\CH^2(\bA)_{\alg,\,\ltor} @=          \CH^2(\bA)_{\alg,\,\ltor}
\end{CD}
\end{equation}
 where the first vertical arrow is the natural inclusion.  Then
Theorem~\ref{merkurev-suslin} implies that $\lambda_2$ is injective on
$\ell$-torsion, from which it follows that $\lambda'_2$ must also be
injective.  On the other hand, Theorem~\ref{suwa} tells us that
$\lambda'_2$ is surjective.  It follows that $\lambda'_2$, and indeed,
all the arrows in diagram~\eqref{diagram}, must be isomorphisms. 
Thus the two groups in \eqref{torgps} are equal, and the Griffiths group
$\Gr^2(\bA)$ has no $\ell$-primary torsion for any $\ell\neq p$. This
completes the proof.
\end{proof}

\begin{corollary}
We have in the notation of Theorem~\ref{main}{\rm :}
\begin{equation}
H^3_{et}(A_{\bk},\Q_{\ell}/\Z_\ell(2)) \isom
\CH^2_{\hom}(A_{\bk})\tensor\Z_\ell  \isom 
N^1H^3_{et}(A_{\bk},\Q_{\ell}/\Z_\ell(2)),
\end{equation}
where $N^1$ denotes the first step of the coniveau filtration.
\end{corollary}
\begin{proof}
The first isomorphism follows from the proof of \ref{main}, while the
second can be found in \cite[section 18]{merkurev83} (also see
\cite{ctr}).
\end{proof}

	The method of proof of Theorem~\ref{main} also proves the
corresponding result for supersingular Fermat threefolds over a finite
field $k$. Recall that a smooth, projective Fermat threefold
$X\subset\P^4$ is said to be supersingular if $H^3_{\rm cris}(X/W(k))$
is of slope $3/2$ (see \cite[Th\'eor\`eme 4.8.1, page 88]{suwa88},
\cite{gouvea-yui}).

\begin{thm} \label{fermat}
	Let $k$ be a finite field of characteristic $p$ and let $\bk$
be an algebraic closure of~$k$.  Let $X$ be a supersingular Fermat
threefold over~$k$ of degree~$m$. Then the Griffiths group of
codimension two cycles on $X_{\bk}$ is at most $p$-primary torsion.
\end{thm}
\begin{proof}
The diagram~\eqref{diagram} is also valid for a smooth Fermat
threefold. As we are over a finite ground field the result of
\cite[Th\'eor\`eme 3]{soule84} applies and so $\CH^2(X_{\bk})_{\hom}$
is torsion and by \cite[Th\'eor\`eme 4.8.1, page 88]{suwa88} the map
\[
 \lambda_2:\CH^2(X_\bk)\to H^3_{et}(X_{\bk},\Q_\ell/\Z_\ell(2))
\]
 is surjective as $X$ is supersingular.  Then we are done by
\ref{merkurev-suslin}.
\end{proof}

\section{Ordinary Abelian Threefolds}
	In order to provide a contrast to the results in the
supersingular case, in this section we observe that the behaviour of
the Griffiths group is controlled by the slope filtration.  This idea
goes back to Bloch (see \cite[Lecture 6, page
6.12]{bloch-lectures}). Our remarks, which are no doubt well-known to
experts, are inspired by the work of Schoen (see \cite{schoen99}).
Recall that an abelian variety over a field $k$ is said to be ordinary
if its Hodge and Newton polygons (for the first crystalline
cohomology) coincide. Equivalently an abelian variety $A$ is ordinary
if and only if the group of geometric points of order $p$ has order
$p^{\dim(A)}$.

Recall that for any smooth, projective variety $X$ over $k$ there are
Abel-Jacobi maps (see \cite{jannsen-mixed}, \cite{schoen95})
\begin{equation}\label{abel-jacobi}
\alpha^r: \CH^r(\bX) \to J_\ell^r(X),
\end{equation}
where 
\[
J_\ell^r(X):=\lim_{k'/k} H^1(G_{k'},
H^{2r-1}(\bX,\Z_\ell(r))/{\rm Tors}),
\]
  the limit is taken over finite Galois extensions~$k'$ of~$k$, and
the cohomology is the continuous Galois cohomology (i.e. cocycles are
continuous with respect to the topology on the Galois group and the
$\ell$-adic topology on the Galois modules). When $k$ is a finite
field $J_\ell^r(X)$ is a torsion group (see \cite[Lemma 1.4,page
4]{schoen95}).

	Before we begin, we remind the reader of the following variant
of Bloch's result (see \cite[Lecture 1]{bloch-lectures}).  This result
is implicit in \cite[Lecture 1]{bloch-lectures}---we give a proof here
for completeness (as we don't know any explicit reference) as it
indicates the relation between the slope filtration and the behaviour
of the Chow groups.
\begin{thm}
Let $X/k$ be a smooth, projective surface over an uncountable
algebraically closed field of characteristic $p$. Further assume that
$H^2_{\rm cris}(X/W(k))$ has a non-trivial slope zero part. Then
$CH^2(X_k)$ is not representable. 
\end{thm}
\begin{proof}
	By \cite[Lecture 1]{bloch-lectures}, it suffices to prove that
the hypothesis imply that the group of transcendental cycles in
\'etale cohomology is non-trivial.  Assume, if possible that it is
trivial, that is, the group $H^2_{\rm
et}(X,\Q_\ell(1))/\text{image}(NS(X_k))_{\Q_\ell}=0$.  This says that
the cycle class map is surjective.  We may assume that $X$ and a basis
for $NS(X_k)$ are defined over a finitely generated subfield of
$k$. Then by further spreading out to a finitely generated ring as our
base. We can, by shrinking the base if neccessary, assume that all the
fibres are smooth.  Then see that there is an non-empty zariski open
set on the base where the Newton polygon of the second crystalline
cohomology of every special fibre coincides with the Newton polygon of
the generic fibre and hence has a non-trivial slope zero part (this
follows from a theorem of Katz and Grothendieck \cite{katz79}). For
any such special fibre, which is defined over a finite field. By
\cite{katz74} we know that over a finite field the characteristic
polynomial of frobenius on the $\ell$-adic cohomology coincides with
the characteristic polynomial of Frobenius on crystalline
cohomology. Hence the cycle class map from the Neron-Severi group to
the second crystalline cohomology cannot be surjective as the
crystalline cohomology has a non-trivial slope zero part (as the image
of Neron-Severi group is contained in the slope 1 part of the
crystalline cohomology). Thus one has a contradiction as the rank of
Neron-Severi does not decrease under specialization.
\end{proof}

\begin{proposition} Let $k$ be a finite field of characteristic $p>2$.
Assume that the Tate conjecture is valid for all smooth projective
surfaces and for all finite extensions of~$k$.  Then the Griffiths
group of any smooth projective, ordinary abelian threefold $A$ over
$\bk$ is non-trivial. More precisely, for all but finite number of
primes $\ell\neq p$, $\Gr^2(\bA)\tensor\Z_\ell\neq 0$.
\end{proposition} 

\begin{proof}
 Suppose $A/k$ is an ordinary abelian threefold.  From the work of
Soul\'e (see \cite[Th\'eor\`eme 3]{soule84}) we know that the Chow
group of homologically trivial cycles on an abelian threefold is
torsion.  Thus the Griffiths group is torsion as well.  Then by
\cite[Corollaire 3.4, page 83]{suwa88}, as $\bA$ is
ordinary the map
\begin{equation}
\lambda_2:\CH^2(\bA)_{\alg,\ltor} 
\to H^3_{et}(\bA,\Q_\ell/\Z_\ell(2))
\end{equation}
 cannot be surjective because it has a non-trivial slope zero part in
$H_{\rm crys}^3(A_{\bk}/W(\bk))$.

 On the other hand we know from the work of Schoen (see \cite[Theorem
0.1, page 795]{schoen99}) that the Abel-Jacobi map $\alpha^2$ is
surjective, for all but finite number of primes $\ell\neq p$, under
the assumption that the Tate conjecture holds for all smooth
projective surfaces.  Thus it suffices to verify that the maps
$\lambda_2$ and $\alpha^2$ coincide on homologically trivial
$\ell$-power torsion cycles, which in turn follows from the
construction of the map $\lambda_2$ given by Raskind (see
\cite[Section 2]{raskind86}).
\end{proof}

\begin{remark}
We would like to complement the above proposition with the following
example which illustrates that the $p$-torsion in the Griffiths group
may be zero even when the abelian variety is an ordinary abelian
variety over a finite field $k$ of characteristic $p>0$. Let $E/k$ be
an ordinary elliptic curve. It is standard result of Deuring that $E$
admits a lifting to an elliptic curve $C$ with complex multiplication
defined over a number field (for a modern proof see \cite[page
192]{oort87}). Let $A=E\times_k E\times_k E$. Then $A$ is an ordinary
abelian variety. The entire discussion in \cite[Remark~4.3, page
601]{gros88} goes through for $A$, and one has that the $p$-adic
Abel-Jacobi mapping constructed
\begin{equation}
\CH^2(A_{\bk})_{\alg}\tensor\Z_p\to H^3_{\rm log}(A,\Q_p/\Z_p(2))
\end{equation}
where the target is the logarithmic cohomology (see \cite{gros88} for
the notation and terminology) is surjective. By \cite[Th\'eor\`eme
3]{soule84} we know that the the kernel of the crystalline cycle class
map is torsion as $X$ is an abelian threefold over a finite
field. Hence we can apply the argument given above to deduce that the
$\Gr^2(A_k)\tensor\Z_p=0$. Thus $p$-torsion homologically trivial
cycles may fail to carry a non-trivial filtration even in presence of
non-trivial slope filtration.
\end{remark}

\bibliographystyle{plain}

\begin{thebibliography}{10}

\bibitem{beauville86}
A.~Beauville.
\newblock Sur l'anneau de {C}how d'une vari\'et\'e ab\'elienne.
\newblock {\em Math. {Ann}.}, 273(4):647--651, 1986.

\bibitem{bloch79}
S.~Bloch.
\newblock Torsion algebraic cycles and a theorem of {R}oitman.
\newblock {\em Compositio {M}ath.}, 39:107--127, 1979.

\bibitem{bloch-lectures}
S.~Bloch.
\newblock {\em Lectures on algebraic cycles}.
\newblock Duke {University} {Math}. {series}. Duke {U}niversity {P}ress,
1980.

\bibitem{colliot-thelene83} 
J.-L.~Colliot-Th\'el\`ene and J.-J.~Sansuc and Ch. Soul\'e.
\newblock Torsion dans le groupe de Chow de codimension deux.
\newblock {\em Duke {Math}. Journal}, 50(3):763--801, 1983.


\bibitem{ctr}
J.-L.~Colliot-Th\'el\`ene and W.~Raskind.
\newblock Groupe de {C}how de codimension deux des vari\'et\`es 
d\'efinies sur
  un corps de nombres: un th\'eor\`eme de finitude pour la torsion.
\newblock {\em Invent. {Math}.}, 105:221--245, 1991.

\bibitem{deninger91}
C.~Deninger and {J}. Murre.
\newblock Motivic decomposition of abelian schemes and the {F}ourier
transform.
\newblock {\em J.~Reine {A}ngew. {Math}.}, 422:201--219, 1991.

\bibitem{gros88}
M.~Gros and {N}.~Suwa.
\newblock Application d'{Abel}-{Jacobi} $p$-adique et cycles
alg\'ebriques.
\newblock {\em Duke {Math}. Journal}, 57(2):578-613.

\bibitem{gouvea-yui} F.~Gouvea and {N}.~Yui.  
\newblock Arithmetic of diagonal hypersurfaces over finite fields.  
\newblock {\em London {M}ath. {S}oc. {L}ecture {N}otes 209}, 
{C}ambridge {U}niversity {P}ress, {C}ambridge, 1995.


\bibitem{fakhruddin95}
N.~Fakhruddin.
\newblock Remarks on the {Chow} groups of supersingular varieties.
\newblock Preprint, 1995.

\bibitem{jannsen-mixed}
U.~Jannsen.
\newblock {\em Mixed motives and algebraic $K$-theory}, volume 1400 of
 {\em  Lecture {N}otes in {Math}.}
\newblock Springer-{V}erlag, {B}erlin, 1990.

\bibitem{katz74}
N.~Katz.
\newblock Some consequences of the {R}iemann hypothesis for 
varieties over finite fields.
\newblock {\em Invent. {M}ath.}, 23:73-77, 1974.

\bibitem{katz79}
N.~Katz.
\newblock Slope filtration of $F$-crystals.
\newblock {\em Ast\'erisque}, 63:113-164, 1979.

\bibitem{kunneman93}
K.~K{\"u}nneman.
\newblock A {L}efschetz decomposition for {C}how motives of abelian
schemes.
\newblock {\em Invent. {Math}.}, 113(1):85--102, 1993.

\bibitem{merkurev83}
A.~S. Merkur'ev and A.~A. Suslin.
\newblock $k$-cohmology of {S}everi-{B}rauer varieties and the norm
residue   homomorphism.
\newblock {\em Math. {U.S.S.R} {I}zvestiya}, 21:307--340, 1983.

\bibitem{mukai81}
S.~Mukai.
\newblock Duality between ${D}({X})$ and ${D}(\hat{X})$ with its
applications to  {P}icard sheaves.
\newblock {\em {N}agoya {M}ath. Journal}, 81:153--175, 1981.

\bibitem{oort74}
F.~Oort.
\newblock Subvarieties of moduli spaces.
\newblock {\em {I}nvent. {M}ath.}, 24:95-119, 1974.

\bibitem{oort87}
F.~Oort.
\newblock Lifting algebraice curves, abelian varieties and their
endomorphisms. 
\newblock {\em {A}lgebraic {G}eometry}, Proceedings {S}ymp. {P}ure
{M}ath., 46 Vol {II}, 165-195, 1987.
 
\bibitem{raskind86}
W.~{R}askind.
\newblock A finiteness theorem in the {G}alois cohomology of algebraic
number  fields.
\newblock {\em Trans. {A}mer. {M}ath. {S}oc.}, 63:107--152, 1986.

\bibitem{schoen95}
C.~Schoen.
\newblock On the computation of the cycle class map for nullhomologous
cycles over the algebraic closure of a finite field.
\newblock {\em Ann. {Scient}. {\'Ecole} {Norm}. {Sup}. $4^{er}$ s\'erie},
  28:1--50, 1995.

\bibitem{schoen99}
C.~Schoen.
\newblock On the image of the $\ell$-adic {A}bel-{J}acobi map for
a variety over the algebraic closure of a finite field.
\newblock {\em J. {A}mer. {M}ath. {Soc}.}, 12(3):795--838, 1999.

\bibitem{soule84}
C.~Soul{\'e}.
\newblock Groupes de {C}how et {K}-th\'eorie de vari\'et\'es sur
un corps fini.
\newblock {\em Math. {A}nn.}, 268:317--345, 1984.

\bibitem{suwa88}
N.~Suwa.
\newblock Sur l'image de l'application d'{A}bel-{Jacobi} de {B}loch.
\newblock {\em Bull. {Soc}. {Math}. {France}}, 116:69--101, 1988.

\bibitem{tate65}
J.~Tate.
\newblock Algebraic cycles and poles of zeta functions.
\newblock In {\em Arithmetical {A}lgebraic {G}eometry ({Proc. Conf.
Purdue.  Univ})}, pages 93--110, New {Y}ork, 1965. Purdue {U}niv.,
Harper \& {Row}.

\end{thebibliography}

\end{document}